\documentclass[aos,MSNbibl,nameyear,seceqn,dvips]{arximspdf}
\usepackage{graphicx}
\doi{10.1214/13-AOS1117} 
\volume{41}
\issue{1}
\pubyear{2013}
\firstpage{1693}
\lastpage{1715}

\makeatletter

\newcommand{\cal}{\mathcal}
\def\real{\mathbb{R}}

\newcommand{\ds}{\displaystyle}
\newcommand{\CCp}{\mathcal{P}}
\newcommand{\CCs}{\mathcal{S}}

\newtheorem{theorem}{Theorem}[section]

\newproclaim{example}{Example}[section]
\newproclaim{remark}{Remark}[section]

\newcommand{\X}{{{\cal X}}}
\newcommand{\B}{{{\cal B}}}

\makeatother

\begin{document}
\begin{frontmatter}

\title{Robust $T$-optimal discriminating designs\thanksref{T1}}
\runtitle{Robust discriminating designs}

\thankstext{T1}{Supported in part by the Collaborative Research Center
``Statistical modeling of nonlinear dynamic processes'' (SFB 823,
Teilprojekt C2) of the German Research Foundation (DFG).}

\begin{aug}
\author[A]{\fnms{Holger} \snm{Dette}\corref{}\ead[label=e1]{holger.dette@rub.de}},
\author[B]{\fnms{Viatcheslav B.} \snm{Melas}\thanksref{t2}\ead[label=e2]{vbmelas@post.ru}}
\and
\author[B]{\fnms{Petr} \snm{Shpilev}\thanksref{t2}}
\runauthor{H. Dette, V. B. Melas and P. Shpilev}
\affiliation{Ruhr-Universit\"at Bochum, St. Petersburg State University
and St.~Petersburg State University}
\address[A]{H. Dette\\
Fakult\"{a}t f\"{u}r Mathematik\\
Ruhr-Universit\"{a}t Bochum\\
44780 Bochum\\
Germany\\
\printead{e1}}
\address[B]{V. B. Melas\\
P. Shpilev\\
Faculty of Mathematics \& Mechanics\\
St. Petersburg State University\\
28, Universitetsky Avenue, Petrodvorets\\
198504 St. Petersburg\\
Russia\\
\printead{e2}}
\end{aug}

\thankstext{t2}{Supported by Russian
Foundation for Basic Research (project number 12-01-00747).}

\received{\smonth{7} \syear{2012}}
\revised{\smonth{1} \syear{2013}}

%
\begin{abstract}
This paper considers the problem of constructing optimal discriminating
experimental designs for competing regression models on the basis of
the $T$-optimality criterion introduced by Atkinson and Fedorov
[\textit{Biometrika} \textbf{62} (\citeyear{atkfed1975a}) 57--70].
$T$-optimal designs depend on unknown model parameters and it is
demonstrated that these designs are sensitive with respect to
misspecification. As a solution to this problem we propose a Bayesian
and standardized maximin approach to construct robust and efficient
discriminating designs on the basis of the $T$-optimality criterion. It
is shown that the corresponding Bayesian and standardized maximin
optimality criteria are closely related to linear optimality
criteria. For the problem of discriminating between two polynomial
regression models which differ in the degree by two the robust
$T$-optimal discriminating designs can be found explicitly. The
results are illustrated in several examples.
\end{abstract}

%

\end{frontmatter}

\section{Introduction} \label{sec1}

An important problem of regression analysis is the identification
of an appropriate model to describe the relation between the response
and a
predictor. Typical examples include dose response studies [see, e.g.,
\citet{brepinbra2005}] in medicine or toxicology or problems in
pharmacokinetics, where a model has usually to be chosen from a class of
competing regression functions; see, for example,
\citet{atkbogbog1998}, \citet{aspmac2000}, \citet{ucibog2005} or \citet
{fooduf2011}. Because a misspecification of a
regression model can result in an inefficient---in the worst case,
incorrect---data analysis, several authors argue that the design of the
experiment
should take the problem of model identification into account. Meanwhile
a huge amount of literature
can be found which addresses the construction of efficient designs for
model discrimination. The literature can be roughly decomposed into two groups.

\citet{hunrei1965}, \citet{stigler1971}, \citet{hill1978}, \citet
{studden1982}, \citet{spruill1990}, Dette (\citeyear{dette1994a,dette1995b}),
\citet{dethal1998},
\citet{songwong1999}, \citet{detmelwong2005} (among many others)
considered two nested models, where
the extended model reduces to the ``smaller'' model for a specific
choice of a subset of the
parameters. The optimal discriminating designs are then constructed
such that
these parameters are estimated most precisely. This concept relies heavily
on the assumption of nested models, and as an alternative
\citet{atkfed1975a}
introduced in a fundamental paper the $T$-optimality criterion for
discriminating between two competing
regression models. Since its introduction this criterion has been
studied by numerous authors
[\citet{atkfed1975b}, \citet{ponce1991} \citet{ucibog2005}, \citet{wawoecle08},
\citet{dettit2009}, \citet{atkinson2010}, \citet{tomlop2010},
Wiens (\citeyear{wiens2009,wiens2010})
or \citet{detmelshp2012} among others].

The $T$-optimal design problem is essentially a maximin problem, and
the criterion can also be applied for nonnested models. Except for very
simple models, $T$-optimal discriminating designs are not easy to find
and even their numerical determination is a very challenging task.
Moreover, an important drawback of this approach consists of the fact
that the criterion and, as a consequence, the corresponding optimal
discriminating designs depend sensitively on the parameters of one of
the competing regression models. In contrast to other optimality
criteria this dependence appears even in the case where only linear
models have to be discriminated. Therefore $T$-optimal designs are
locally optimal in the sense of \citet{chernoff1953} as they can only be
implemented if some prior information regarding these parameters is
available. Moreover, we will demonstrate in Example~\ref{emax} that the
efficiency of a $T$-optimal design depends sensitively on a precise
specification of the unknown parameters in the criterion. This problem
has already been recognized by \citet{atkfed1975a} who proposed Bayesian
or minimax versions of the $T$-optimality criterion. However---to the
best knowledge of the authors---there exist no results in the
literature investigating optimal design problems of this type more
rigorously (we are not even aware of any numerical solutions).

The present paper is devoted to a more detailed discussion
of robust $T$-optimal discriminating designs.
We will study a Bayesian and a standardized maximin
version of the $T$-optimal discriminating design problem; see \citet
{chaver1995} and \citet{dette1997}. It is demonstrated that optimal
designs with respect to these criteria are closely related to optimal
designs with respect to linear optimality criteria. For the particular
case of discriminating between two competing polynomial regression models
which differ in the degree by two, robust $T$-optimal discriminating designs
are found explicitly. These results provide---to our best of our
knowledge---the first explicit solution
in this context. Interestingly, the structure of these Bayesian and
standardized maximin $T$-optimal discriminating
designs is closely related to the structure of designs for a most
precise estimation of
the two highest coefficients
in a polynomial regression model; see \citet{gaffke1987} or \citet{studden1989}.

The remaining part of the paper is organized as follows. In Section
\ref{sec2} we revisit the $T$-optimality criterion introduced by
\citet{atkfed1975a} for two regression models, which will be called
locally $T$-optimality criterion in order to reflect the dependence on
the parameters of one of the competing models. In particular it is
demonstrated that locally $T$-optimal designs can be inefficient if the
parameters in the optimality criterion have been misspecified.
Section~\ref{sec3} is devoted to robust versions of the $T$-optimality
criterion and properties of the corresponding optimal designs, while
Section~\ref{sec4} gives explicit results for Bayesian and standardized
maximin $T$-optimal discriminating designs for two competing polynomial
regression models. In Section~\ref{sec5} we illustrate the results and
construct robust optimal discriminating designs for a constant and
quadratic regression. These two models have been proposed in
\citet{brepinbra2005} to detect dose response signal in phase II
clinical trial if there is some evidence that the shape of the dose
response might be u-shaped. In particular it is demonstrated by a small
simulation study that the $T$-optimal discriminating designs improve
the power of the $F$-test for discriminating between the two polynomial
models.

\section{Locally $T$-optimal designs} \label{sec2}

We assume that the relation between a predictor $x$ and response $y$
is described by the regression model
\[
y=\eta(x)+\varepsilon,
\]
where $ x$ varies in a compact designs space ${\X}\subset\mathbb
{R}^k$, and $\varepsilon$
denotes a centered random variable with finite variance. We also assume
that observations at experimental conditions $x_1$ and $x_2$ are
independent and that
there exist two competing continuous parametric models,
say $\eta_1 $ or $\eta_2$, for the regression function $\eta$ with
corresponding parameters
$\theta_1\in\mathbb{R}^{m_1}$; $\theta_2\in
\mathbb{R}^{m_2}$, respectively.
In order to find ``good'' designs for discriminating between the models
$\eta_1$ and $\eta_2$, we consider approximate designs in the
sense of \citet{Kiefer1974},
which are defined as probability measures on the design space $\mathcal
{X}$ with finite support.
The support points, say $x_1,\ldots, x_s$, of an (approximate) design
$\xi$ give the locations where observations
are taken, while the weights give the corresponding relative
proportions of total observations to be taken at these points.
If the design $\xi$ has masses $\omega_i>0 $ at the different points
$x_i$ $(i =
1, \ldots, s)$ and $n$ observations can be made by the experimenter,
the quantities
$\omega_i n$ are rounded to integers, say $n_i$, satisfying $\sum
^s_{i=1} n_i =n$, and
the experimenter takes $n_i$ observations at each location $x_i$
$(i=1, \ldots, s)$.
It has been demonstrated by \citet{pukrie1992} that the loss of
efficiency caused by rounding is of order $O(n^{-2})$ and $O(n^{-1})$
for differentiable and nondifferentiable optimality criteria,
respectively.

To determine a good design for discriminating between the two rival
regression models $\eta_1$ and $\eta_2$, \citet{atkfed1975a}
proposed in a fundamental paper to fix one model, say $\eta_2$ (more
precisely its corresponding parameter $\theta_2$), and to determine
the design which maximizes the minimal deviation
%
\begin{equation}
\label{topt} T(\xi,\theta_2) = \min_{\theta_1 \in\Theta_1 }
\int_{\chi
} \bigl(\eta_1(x,\theta_1)-
\eta_2 (x,\theta_2) \bigr)^{2}\xi(d x)
\end{equation}
between the model $\eta_2$ and the class of models $\{\eta_1(x,\theta
_1) |\theta_1 \in\Theta_1 \}$ defined by $\eta_1$, that is,
%
\begin{equation}
\label{tloc} \xi^*=\mathop{\arg\max}_{\xi} T(\xi,\theta_2).
\end{equation}
Note that a design maximizing (\ref{topt}) is constructed, such that
the deviation between the given model $\eta_2(\cdot,\theta_2)$ and
its best approximation by models of the form $\eta_1(\cdot,\theta
_1)$ with respect to the $L^2(\xi)$-distance is maximal. Moreover, in
linear and nested models $\eta_1$ and $\eta_2$ it can be shown that
the design $\xi^*$ in (\ref{teff}) maximizes the power of the
corresponding $F$-test. For these and further properties of the
criterion we refer to \citet{dettit2009}.
Throughout this paper we call the maximizing design and optimality
criterion in (\ref{tloc}) locally $T$-optimal discriminating design
and local
$T$-optimality criterion, respectively, because they will depend on the
specification of the parameter $\theta_2$ used for the model
$\eta_2$. The local $T$-optimal design problem is a maximin problem
and except for very
simple models the corresponding optimal designs are extremely hard to
find. Even their numerical construction is a difficult and challenging task.
Nevertheless, since its introduction the optimal designs with respect
to the criterion
(\ref{topt}) have found considerable interest in the literature and we refer
the interested reader to the work of \citet{ucibog2005} or \citet
{dettit2009} among others. The latter authors showed
that the optimization problem
(\ref{tloc}) is closely related to a problem in nonlinear
approximation theory, that is,
%
\begin{equation}
\label{rdef} R(\theta_2):= \max_\xi T(\xi,
\theta_2) = \inf_{\theta_1 \in
\Theta_1} \sup_{x \in\mathcal{X}}
\bigl|\eta_1(x,\theta_1)-\eta_2(x,
\theta_2)\bigr|^2,
\end{equation}
where $T(\xi,\theta_2)$ is defined in (\ref{topt}).
Because of its local character, locally $T$-optimal designs are rather
sensitive with respect to the misspecification
of the unknown parameter and the following example illustrates this fact.

\begin{example} \label{emax}
We consider the problem of constructing a $T$-optimal discriminating
design for the Michaelis--Menten model,
\[
\eta_1 (x,\theta_1)=\frac{\theta_{1,1}x}{\theta_{1,2}+x}
\]
[see, e.g., \citet{cornish1995}] and the EMAX model
\[
\eta_2 (x,\theta_2)=\theta_{2,0}+
\frac{\theta_{2,1}x}{\theta_{2,2}+x};
\]
see, for example, \citet{Danesi2002}. It is easy to see that the
$T$-optimal \mbox{discriminating} design does not depend on the parameter
$\theta_{2,1}$ and therefore we assume without loss of generality
$\theta_{2,1}\equiv1$. In Table~\ref{table1} we display some locally
%
\begin{table}
\tabcolsep=0pt
\caption{The support points and the weights of $T$-optimal
discriminating designs for a Michaelis--Menten and an EMAX model and
various specifications of the parameters $\theta_{2,0}$ and
$\theta_{2,2}$ of the EMAX model. The locally $T$-optimal design puts
weights $w_1, w_2$ and $w_3$ at the points $1$, $x^*$ and $2$,
respectively} \label{table1}
\begin{tabular*}{\tablewidth}{@{\extracolsep{\fill}}lccccccccccc@{}}
\hline
$\bolds{\theta_{2,0}}$ & $\bolds{\theta_{2,2}}$ & $\bolds{x^*}$
& $\bolds{w_1}$ & $\bolds{w_2}$ & $\bolds{w_3}$ &
$\bolds{\theta_{2,0}}$ & $\bolds{\theta_{2,2}}$
& $\bolds{x^*}$ & $\bolds{w_1}$ & $\bolds{w_2}$ & $\bolds{w_3}$ \\
\hline
$-2$ & $2$ & $1.368$ & $0.206$ & $0.499$ & $0.295$ & $-2$ & $1$ &
$1.352$ & $0.211$ & $0.499$ & $0.29$\hphantom{0} \\
$-1$ & $2$ & $1.347$ & $0.176$ & $0.495$ & $0.329$ & $-1$ & $1$ &
$1.321$ & $0.165$ & $0.491$ & $0.344$ \\
$-1/2$ & $2$ & $1.211$ & $0.040$ & $0.584$ & $0.376$ & $-1/2$ & $1$ &
$1.590$ & $0.619$ & $0.336$ & $0.045$\\
$1/2$ & $2$ & $1.400$ & $0.260$ & $0.498$ & $0.242$ & $1/2$ & $1$ &
$1.384$ & $0.261$ & $0.498$ & $0.239$ \\
$1$ & $2$ & $1.390$ & $0.247$ & $0.499$ & $0.254$ & $1$ & $1$ & $1.378$
& $0.253$ & $0.499$ & $0.248$ \\
$2$ & $2$ & $1.387$ & $0.238$ & $0.499$ & $0.263$ & $2$ & $1$ & $1.337$
& $0.244$ & $0.500$ & $0.256$ \\
\hline
\end{tabular*}
\end{table}
$T$-optimal discriminating designs on the interval $[1,2]$ for various
values of parameters $\theta_{2,i}, i=0,2$. We observe that the
resulting designs are rather sensitive with respect to the
specification of the values $\theta_{2,0}$ and $\theta_{2,2}$. Note
that in contrast to the $T$-optimal discriminating design, the
$T$-efficiency
%
\begin{equation}
\label{teff} \mathrm{Eff}_T (\xi,\theta_2) =
\frac{T(\xi,\theta_2)}{\sup_\eta
T(\eta,\theta_2)}
\end{equation}
depends also on the parameter $\theta_{2,1}$ of the EMAX model and
some efficiencies are depicted in Figure~\ref{Pic4}
if the true values are given by $\theta_{2,0}=-1$, $\theta_{2,1}=1$,
$\theta_{2,2} \in(2,6)$, and
one uses the $T$-optimal discriminating design calculated under the
%
\begin{figure}

\includegraphics{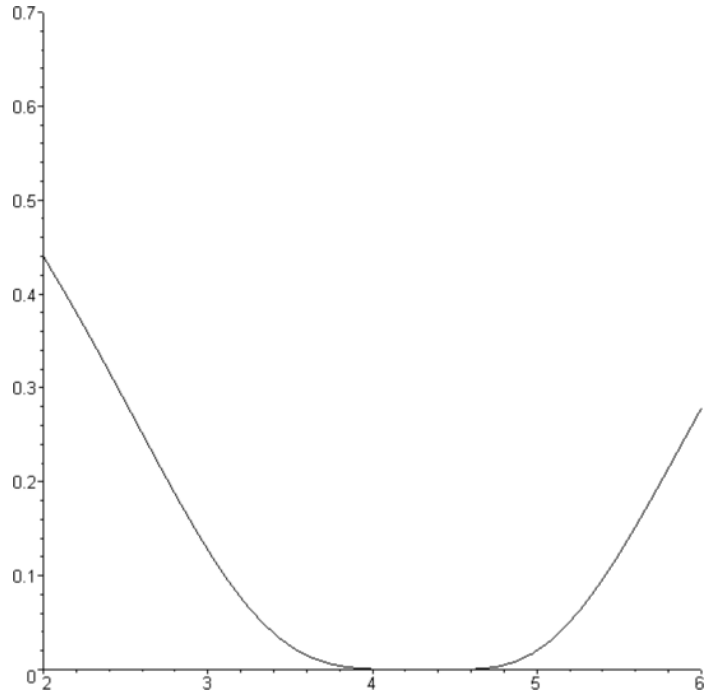}

\caption{$T$-efficiency (\protect\ref{teff}) of the locally
$T$-optimal discriminating design for Michaelis--Menten and EMAX
model calculated under the assumption $\theta_{2,0}=-1/4$,
$\theta_{2,1}=1$ while the ``true'' values are given by
$\theta_{2,0}=-1$, $\theta_{2,1}=1$. The efficiencies depend on
the parameter $\theta_{2,2} \in(2,6) $.}
\label{Pic4}
\end{figure}
assumption $\theta_{2,0}=-1/4$, $\theta_{2,1}=1$ and
{$\theta_{2,2}\in(2,6) $}. We observe
a substantial loss of $T$-efficiency in some regions for
$\theta_{2,2}$. If $\theta_{2,2} \in(0,2) $ the efficiency is
larger than $50\%$, if $\theta_{2,2} \in(2,3) \cup(5.5,6)$ it varies
between $15 \%$ and $40 \% $, if $\theta_{2,2} \in(3, 5.5)$
the efficiency is smaller than $15\%$, and the locally $T$-optimal
design cannot be recommended. On the basis of these observations it
might be desirable to use
designs which are less sensitive with respect to misspecification of
the parameter $\theta_{2,0}$, and the
corresponding methodology will be developed\vadjust{\goodbreak} in the following
section. Robust $T$-optimal designs for discriminating between the
Michaelis and EMAX model will be discussed at the end of this paper
where we construct a uniformly better design; see Section~\ref{sec53}.
\end{example}

\section{Robust $T$-optimal discriminating designs} \label{sec3}

Because the previous example indicates that locally $T$-optimal
discriminating designs are
sensitive with respect to misspecification of the parameters $\theta
_2$ of the model $\eta_2$ in the $T$-optimality criterion (\ref{topt}),
the consideration of robust optimality criteria for model
discrimination is
of great interest. In the context of constructing efficient robust
designs for parameter
estimation in nonlinear regression models Bayesian and standardized
maximin optimality
criteria have been discussed intensively in the literature; see \citet
{chaver1995}, \citet{dette1997} or \citet{muepaz1998},
among many others. However, to our best knowledge, these methods have
not been investigated rigorously in the context of
model discrimination so far, and in this
section we will define a robust version of the local $T$-optimality criterion.
Recall the definition of this criterion in (\ref{topt}) and its
optimal value $R(\theta_2)$ in (\ref{rdef}); then a design $\xi^*_M$
is called standardized maximin $T$-optimal discriminating (with respect
to the set $\Theta_2$) if it maximizes the criterion
%
\begin{equation}
\label{critbay} V_M (\xi) = \inf_{\theta_2 \in\Theta_2}
\frac{T(\xi,\theta
_2)}{R(\theta_2)},
\end{equation}
where $\Theta_2$ is a pre-specified set, reflecting the experimenter's
belief about the unknown parameter $\theta_2$. Similarly, if $\pi$
denotes a prior distribution on the set $\Theta_2$, then a design $\xi
^*_B$ is called Bayesian $T$-optimal (with respect to the prior $\pi$)
if it maximizes the criterion
%
\begin{equation}
\label{critvb} V_B(\xi) = \int_{\Theta_2}
T(\xi,\theta_2)\pi(d \theta_2).
\end{equation}
We would like to point out here that criteria (\ref{critbay}) and
(\ref{critvb}) yield to different optimal designs. Some interesting
relations between both optimality criteria can be found in \citet
{dethaiimh2007}. In particular, these authors showed that, under
appropriate regularity assumptions, a design maximizing criterion
(\ref{critbay}) is always a Bayesian $T$-optimal design with respect
to a prior supported at the extreme points of the best approximation of
the function $\eta_2(\cdot, \theta_2)$ by the function $\eta
_1(\cdot, \theta_1)$ with respect to the sup-norm $\| g \|_\infty=
\sup_{x \in\mathcal{X}}| g(x)|$. From a practical point of view the
use of criterion (\ref{critbay}) or (\ref{critvb}) is a matter of
taste and depends on the concrete application. If some regions of the
parameter space are more likely than others a Bayesian criterion could
be preferred; otherwise the maximin or Bayesian criterion with an
uninformative prior might be appropriate.

In the following discussion we investigate the problem of constructing
robust discriminating designs for two linear
regression models
%
\begin{equation}
\label{linmod} \eta_1(x,\theta_1)=\sum
^{m_1-1}_{i=0}\theta_{1,i}f_i(x),\qquad
\eta_2(x,\theta_2)=\sum^{m_2-1}_{i=0}
\theta_{2,i}f_i(x),
\end{equation}
where $m_2 > m_1$, $f_0,\ldots,f_{m_2-1}$ are given linearly independent
regression functions, and
$\theta_i=(\theta_{i,0},\ldots,\theta_{i,m_i-1})^T$ denotes the
parameter in the model $\eta_i$ $(i=1,2)$. We introduce the notation
$b_1=\theta_{2,m_1} /\theta_{2,m_2-1},\ldots, b_{m_2-m_1-1} =
\theta_{2,m_2-2} /\theta_{2,m_2-1}$, $m=m_2-1$, $s=m_2-m_1$,
$q_i=\theta_{1,i}-\theta_{2,i}$ ($ i=0,1,\ldots,m-s$) and obtain for
the difference $\eta_1(x,\theta_1) -\eta_2(x,\theta _2)$ the
representation
%
\begin{eqnarray}
\label{help}\qquad
&&
\bar\eta\bigl(x,q,\theta^c_{2,m},
\theta_{2,m}\bigr)\nonumber\\[-8pt]\\[-8pt]
&&\qquad= \sum^{m-s}_{i=0}q_if_i(x)-
\bigl(b_1f_{m-s+1}(x)+\cdots+ b_{s-1}
f_{m-1}(x)+f_m(x)\bigr)\theta_{2,m},\nonumber
\end{eqnarray}
where $\theta^c_{2,m}=(\theta_{2,m-s+1},\ldots,\theta_{2,m-1})^T$.
Thus the locally $T$-optimality criterion in (\ref{topt}) can be
rewritten as
\begin{eqnarray*}
T(\xi, \theta_2) &=& \inf_{\theta_1 \in\mathbb{R}^{m-s+1}} \int
_\mathcal{X} \bigl(\eta_1(x,\theta_1)-
\eta_2(x,\theta_2)\bigr)^2 \,d\xi(x) \\
&=&
\theta^2_{2,m} \inf_{q \in\mathbb{R}^{m-s+1}} \int
_\mathcal{X} \bar\eta^2(x,q,b,1)\,d\xi(x),
\end{eqnarray*}
where $b=(b_1,\ldots,b_{s-1})^T$.
Consequently, locally $T$-optimal designs depend only on the ratios
$b_i=\theta_{2,m-s+i}/\theta_{2,m}$ $(i=1,\ldots,s-1)$. Similarly, if
$\pi$ is a prior distribution for the vector $\theta_2$, then it
follows from these discussions that the Bayesian $T$-optimality
criterion depends only on the induced prior distribution, say $\bar\pi
$, for the parameter $b=(b_1,\ldots,b_{s-1}$).
Therefore we assume that the vector $b$ varies in a subset $\B
\subset\mathbb{R}^{s-1}$ and define $\bar\pi$ as a prior
distribution on $\B$. With these notations the Bayesian
$T$-optimality criterion in (\ref{critvb}) simplifies to
%
\begin{equation}
\label{critbaypol} V_B(\xi) = \int_{\mathcal{B}} \inf
_{q\in
\mathbb{R}^{m-s+1}} \int_{\X} \bar\eta^2(x,q,b,1)
\xi(dx) \bar\pi(db).
\end{equation}
Similarly, we have with the notation $\bar\theta_2=\theta_2/\theta_{2,m}$,
\[
R(\theta_2) = \max_\xi T(\xi,
\theta_2) = \theta^2_{2,m}
R(\bar\theta_2)
\]
and defining
${\cal{B}}=\{(\theta_{2,m-s+1}/\theta_{2,m},\ldots,\theta_{2,m-1}/\theta_{2,m})^T
|\theta_2 \in\Theta_2 \} \subset\mathbb{R}^{s-1}$ and for $b \in\B$
%
\begin{equation}
\label{rquer} \bar R(b)= R\bigl((b_1,\ldots,b_{s-1},1)^T
\bigr)
\end{equation}
the factor $\theta^2_{2,m}$ in (\ref{critbay}) cancels and the
standardized maximin $T$-optimality criterion reduces to
%
\begin{equation}
\label{eq30} V_M(\xi) = \inf_{b\in\B}
\frac{\inf_{q\in\mathbb{R}^{m-s+1}}
\int_{\cal X} \bar\eta^2 (x,q,b,1)\xi(dx)} {
\bar
R(b)} = \inf_{b \in\B} \operatorname{eff}_T (\xi,b),
\end{equation}
where the efficiency is defined in an obvious manner, that is,
\[
\operatorname{eff}_T (\xi,b)=\frac{T(\xi, (b_1,\ldots,b_{s-1},1)^T)}{\bar R(b)}.
\]
Throughout this paper we denote by $f(x)=(f_0(x),
f_1(x),\ldots,f_m(x))^T$ the vector of regression functions with
corresponding decomposition
\begin{eqnarray*}
f_{(1)}(x)&=&\bigl(f_0(x),f_1(x),\ldots,f_{m-s}(x)\bigr)^T \in\mathbb{R}^{m-s+1},
\\
f_{(2)}(x)&=& \bigl(f_{m-s+1}(x),\ldots,f_m(x)
\bigr)^T \in\mathbb{R}^s.
\end{eqnarray*}
We assume that the functions $f_0,\ldots,f_m$ are linearly independent
and continuous on $\mathcal{X}$ and define
\[
M(\xi)=\int_{\mathcal{X}} f(x)f^T(x)\xi(dx)
\]
as the information matrix of a design with corresponding blocks
\[
M_{ij}(\xi)=\int_{\cal X} f_{(i)}(x)f_{(j)}^T(x)
\xi(dx),\qquad i,j=1,2,
\]
and Schur complement
\[
M_{(s)}(\xi)=M_{22}(\xi)-X^TM_{11}(
\xi)X,
\]
where $X \in\mathbb{R}^{m-s+1\times s}$ is an arbitrary solution of
the equation
$M_{11}(\xi)X=M_{12}(\xi)$ [if this equation has no solutions, then
the matrix $M_{(s)}(\xi)$ remains undefined]. Our first main result relates
the Bayesian and standardized maximin $T$-optimality criteria
to linear optimality criteria.
%
\begin{theorem}\label{lem31}
Let $\bar\pi$ denote a prior distribution for the vector $b \in\B$,
such that the
matrix
\[
L=\int_{\mathcal{B}} \pmatrix{ b b^T & b
\cr
b^T & 1} \bar\pi(db)
\]
exists; then the two following statements are equivalent:
\begin{longlist}[(2)]
\item[(1)] The design $\xi^*$ is a Bayesian $T$-optimal
discriminating design with respect to the prior $\bar\pi$ for the
linear regression models defined in (\ref{linmod}).
\item[(2)] The design $\xi^*$ maximizes the linear criterion
\[
\operatorname{tr} LM_{(s)}(\xi)
\]
in the class of all approximate designs $\xi$, for which there exists
a solution $X \in\mathbb{R}^{m-s+1 \times s}$ of the equation
%
\begin{equation}
\label{eq33} M_{11}(\xi)X=M_{12}(\xi).
\end{equation}
\end{longlist}
\end{theorem}

\begin{pf}
If the matrix $M_{(s)}(\xi)$ is nonsingular, then it
follows from Karlin and Studden [(\citeyear{karstu1966}), Section 10.8], that
\[
\bigl(M_{(s)}(\xi) \bigr)^{-1}= \bigl({\mathbf O}^T\dvtx
I_s\bigr)M^{-}(\xi) \pmatrix{{\mathbf O}
\cr
I_s},
\]
where $I_s \in\mathbb{R}^{s \times s}$ is the identity matrix, ${\mathbf
O} \in\mathbb{R}^{m-s+1 \times s}$ is the matrix with all entries
equal to $0$ and $M^{-}(\xi)$ is an arbitrary generalized inverse
of the matrix $M(\xi)$. For any $(m-s+1) \times s$ matrix $K$
we have the inequality
\[
\bigl(-K^T\dvtx  I_s \bigr) \pmatrix{ M_{11}(\xi)&
M_{12}(\xi)
\cr
M_{21}(\xi)&M_{22}(\xi)} \pmatrix{
-K
\cr
I_s} \geq M_{(s)}(\xi),
\]
where there is equality if and only if the matrix $K$ is a solution of
the equation (\ref{eq33}); see \citet{karstu1966},
Section 10.8. From (\ref{help}) and the discussion in the subsequent
paragraph, we obtain the representation
%
\begin{eqnarray}
\label{eq31} T(\xi,\theta_2)&=&\theta^2_{2,m}
\min_{q\in
\mathbb{R}^{m-s+1}}\bigl(q^T,b^T,1\bigr)M(\xi)
\bigl(q^T,b^T,1\bigr)^T \nonumber\\[-8pt]\\[-8pt]
&=&
\theta^2_{2,m} \bigl(b^T,1\bigr)M_{(s)}(
\xi) \bigl(b^T,1\bigr)^T,\nonumber
\end{eqnarray}
where the last equality follows from the fact that each vector
$(q^T,b^T,1)^T$ can be represented in the form
\[
\bigl(q^T,b^T,1\bigr)^T=
\bigl(-K^T\dvtx  I_s\bigr)^T \bigl(b^T,1
\bigr)^T\vadjust{\goodbreak}
\]
for some appropriate matrix $K \in\mathbb{R}^{m-s+1 \times s}$ [just
use the matrix $K=\break-q(b^T,1)/(b^Tb+1)$].
The assertion of Theorem~\ref{lem31} is now obvious.
\end{pf}

A similar result for standardized maximin $T$-optimal discriminating
designs is formulated in the following theorem. Throughout this paper
we will use the notation $\bar{\mathbb{R} } = \mathbb{R}
\cup\{-\infty,\infty\}$ with the usual compactification.

\begin{theorem}\label{lem32}
If $\B\subset\bar{\mathbb{R}}^{s-1} $ be a given compact
set, then the following two statements are equivalent:
\begin{longlist}[(2)]
\item[(1)] The design $\xi^*$ is a standardized maximin $T$-optimal
discriminating design
for the regression models defined in (\ref{linmod}) with respect to
the set $\mathcal{B}$.
\item[(2)] For the design $\xi^*$ there exists a solution to equation
(\ref{eq33}) and a matrix $L^* \in\mathbb{R}^{s \times s}$
such that the pair $(L^*,\xi^*)$ satisfies
%
\begin{eqnarray}
\label{wahlm}
\operatorname{tr} L^*M_{(s)}\bigl(\xi^*\bigr) &=& \sup_{\xi}
\operatorname{tr}L^* M_{(s)}(\xi),
\\
\label{wahl}
\operatorname{tr} L^*M_{(s)}\bigl(\xi^*\bigr) &=& \inf_\mathcal{L}
\operatorname{tr}L M_{(s)}\bigl(\xi^*\bigr),
\end{eqnarray}
where the supremum in (\ref{wahlm}) is taken with respect to all
approximate designs, and the
set ${\cal L}$ in (\ref{wahl}) is defined by
\[
\Biggl\{ \sum^k_{i=1}
\bigl(b_i^T, 1\bigr)^T \bigl(b^T_i,1
\bigr) { \omega_i \over{\bar
R}(b_i)} \Big| b_i\in\B, \omega_i>0,
i=1,\ldots,k,
\sum_{i=1}^k\omega_i = 1
\Biggr\}.
\]
\end{longlist}
\end{theorem}

\begin{pf}
By a similar argument as that used in the proof of Theorem~\ref{lem31},
the standardized $T$-optimality criterion in
(\ref{eq30}) can be represented as
\begin{eqnarray*}
V_M(\xi) &=& \inf_{b \in\mathcal{B}} \inf_{q \in\mathbb
{R}^{m-s+1}}
\frac{(q^T,b^T,1)M(\xi)(q^T,b^T,1)^T}{\bar R(b) }
\\
&=& \inf_{b \in\mathcal{B}} \frac{(b^T,1)M_{(s)}(\xi
)(b^T,1)^T}{\bar R(b)} = \inf_{L \in\mathcal{L}}
\operatorname{tr} LM_{(s)}(\xi).
\end{eqnarray*}
The assertion now follows from the von Neumann theorem on minimax
problems; see \citet{osbrub1994}.
\end{pf}

A lower bound for the efficiencies of a standardized maximin
$T$-optimal discriminating
design is given in the following theorem.

\begin{theorem}\label{theo31}
Let $\xi^*$ denote a standardized maximin $T$-optimal discriminating
design for the linear regression models defined in (\ref{linmod}) with
respect to set $\mathcal{B}$. Then for all $b \in\mathcal{B}$
\[
\operatorname{eff}_T\bigl(\xi^*,b\bigr) \geq\frac{1}{s}.
\]
\end{theorem}

\begin{pf}
Recall the definition of the standardized
maximin optimality criterion in (\ref{eq30}). Because for any $b \in
\mathcal{B}$
\[
\operatorname{eff}_T \bigl(\xi^*,b\bigr) \geq\inf_{b \in\B}
\operatorname{eff}_T \bigl(\xi^*,b\bigr)=V_M\bigl(\xi^*\bigr)
\]
the assertion follows, if the inequality
\[
\sup_{\xi} V_M(\xi)\geq\frac{1}{s}
\]
can be established. For this purpose we define the function
%
\begin{equation}
\label{fctX} \psi(x)=f_{(2)}(x)- X^T f_{(1)}(x),
\end{equation}
where $X$ is an $(m-s+1) \times s$-matrix (the dependence of the
function $\psi$ on this matrix is not reflected in the notation).
Let $\xi$ be an arbitrary design such that the matrix $M_{(s)}(\xi)$
is nonsingular. Then it follows from the Cauchy--Schwarz inequality that
%
\begin{equation}
\label{eq31a}\quad \inf_{l\in\mathbb{R}^s \setminus{0}} \frac
{l^TM_{(s)}(\xi)l}{\sup
_{x\in{\X}}{(l^T\psi(x))^2}} \geq
\frac{1}{\sup_{x\in{\X}}\psi^T(x)M^{-1}_{(s)}(\xi)
\psi(x)}.
\end{equation}
By the equivalence theorem for $D_s$-optimal designs [see
\citet{karstu1966}, Section 10.8] there\vspace*{1pt} exists a design
$\tilde\xi$ and a matrix $\tilde X$ satisfying $M_{11}(\tilde\xi
)\tilde X=M_{12}(\tilde\xi)$,
such that the corresponding matrix $M_{(s)}(\tilde\xi)$ and the
vector $\tilde\psi(x)=f_{(2)}(x) - \tilde X^T f_{(1)}(x)$ satisfy
\[
\max_{x\in{\X}} \tilde\psi^T(x)M^{-1}_{(s)}(
\tilde\xi)\tilde\psi(x) =s.
\]

Consider any design $\xi$ for which a solution $X$ of (\ref{eq33})
exists. Then we have for the corresponding function $\psi$ in (\ref{fctX}),
\[
M_{(s)}(\xi)=\int_{\mathcal{X}} \psi(x)
\psi^T(x)\xi(dx).
\]

Therefore we obtain from formula (\ref{eq31})
\[
\bar R(b) = \theta_{2,m}^2 \max_\xi
\int_{\mathcal{X}} \bigl(\bigl(b^T,1\bigr) \psi(x)
\bigr)^2\xi(dx) = \theta_{2,m}^2\max
_{x\in
\mathcal{X}} \bigl(\bigl(b^T,1\bigr) \psi(x)
\bigr)^2,
\]
which gives for the vector $l=(\theta_{2,m-s+1},\ldots,\theta_{2,m})^T$
%
\[
\sup_{x\in{\X}} \bigl(l^T \psi(x) \bigr)^2=
\bar R(b),
\]
where
$b=(\theta_{2,m-s+1}/\theta_{2,m},\ldots,\theta_{2,m-1}/\theta_{2,m})^T$.
Thus the left-hand side in (\ref{eq31a}) equals $V_M(\xi)$ and
\[
\sup_\xi V_M(\xi)\geq V_M(\tilde
\xi)\geq\frac{1}{s},
\]
which proves the assertion of Theorem~\ref{theo31}.
\end{pf}

\section{Robust $T$-optimal designs for polynomial regression} \label{sec4}

In general locally $T$-optimal discriminating designs have to be found
numerically, and
this statement also applies to the construction of robust $T$-optimal
discriminating designs
with respect to the Bayesian or standardized maximin criterion.
In order to get more insight in the corresponding optimal design
problems we consider in this section the case of two competing
polynomial regression models
which differ in the degree by two. Remarkably, for this situation the
robust $T$-optimal discriminating designs can be found explicitly.
To be precise, let $s=2$, consider the vectors of monomials
\[
f_{(1)}(x)=\bigl(1,x,\ldots,x^{m-2}\bigr)^T,\qquad
f_{(2)}=\bigl(1,x,\ldots,x^{m}\bigr)^T
\]
and define
\[
U_n(x)=\frac{\sin((n+1)\operatorname{arcos} x)}{\sin(\operatorname{arcos} x)}
\]
as the Chebyshev polynomial of the second kind; see \citet{szego1959}.
We assume that the design space is given by the symmetric interval
$[-a,a]$ and consider for $\beta> 0$ designs $\xi_{m,\beta} $
defined as follows. If $ \beta=1$, then the design
$\xi_{m,1}$ puts masses $1/(2(m-1))$ at the points $-a$, $a$ and
masses $1/(m-1)$ at the $m-2$ roots of the polynomial
$U_{m-2}(x/a)$. If $\beta\neq1$ the design $\xi_{m,\beta} $ is
supported at the $m+1$
roots $-a=x_0<x_1<\cdots<x_{m-1}<x_m=a$ of the polynomial
\[
\bigl(x^2-a^2\bigr) \biggl\{U_{m-1} \biggl(
\frac{x}{a} \biggr)+\beta U_{m-3} \biggl(\frac{x}{a}
\biggr) \biggr\},
\]
where the corresponding weights are given by
\begin{eqnarray*}
\xi_{m,\beta}(\mp a)&=&\frac{1+\beta}{2[m+\beta(m-2)]},
\\
\xi_{m,\beta}(x_{j})&=& \biggl[ m-1-\frac{(1+\beta)U_{m-2} (
{x_j}/{a} )}{U_m
({x_j}/{a} )+
\beta
U_{m-2} ({x_j}/{a} )}
\biggr]^{-1},\qquad j=1,\ldots,m-1.
\end{eqnarray*}

\begin{theorem}\label{theo41}
\textup{(1)} Let $\bar\pi$ denote a symmetric prior distribution
on $\B\subseteq(-\infty, \infty)$ with existing second moment, and define
%
\begin{equation}
\label{betab} \beta_B= \min\biggl\{ 1, \frac{\int_\B b^2 \bar\pi(db)}{
a^2} \biggr\}.
\end{equation}
The design $\xi_{m,\beta_B}$ is
a Bayesian $T$-optimal discriminating on the interval $[-a,a]$ for the
polynomial
regression models of degree $m-2$ and $m$.

\textup{(2)} Define $\beta_M=1-2h^*$, where $h^*$ is the unique
maximizer of the function
%
\begin{equation}
\label{eq41} \inf_{b \in\B} \frac{b^2+a^2h}{a^2 \bar R(b,a)} (1-h),
\end{equation}
where
\[
\bar R (b,a)= \inf_{q_0,\ldots,q_{m-2}\in\mathbb{R}} \sup_{x \in
[-1,1]}
a^{2m} \biggl| x^m + \frac{b}{a} x^{m-1} +
q_{m-2} x^{m-2} + \cdots+ q_1 x+q_0
\biggr|^2
\]
in the interval $[0,\frac{1}{2}]$. Then
the design $\xi_{m,\beta_M} $
is a standardized maximin $T$-optimal discriminating
design on the interval $[-a,a]$ for the polynomial
regression models of degree $m-2$ and $m$.
\end{theorem}

\begin{pf}
We will prove the statement using some basic facts of the theory of
canonical moments; see \citet{dettstud1997} for details. To be precise,
let $\CCp([-a,a])$ denote the set of all probability measures on the
interval $[-a,a]$, and denote for a design $\xi\in\CCp([-a,a])$ its
moments by
\[
c_i =c_i(\xi) =\int^a_{-a}
x^i\xi(dx),\qquad i=1,2,\ldots.
\]
Define $\mathcal{M}_k = \{(c_1,\ldots,c_k)^T |\xi\in\CCp ([-a,a])\}$
as the $k$th moment space and $ \Phi_{k}(x) =(x,\ldots,x^{k})$ as the
vector of monomials of order $k$. Consider for a fixed vector $c =(c_1,\ldots, c_{k})^T \in{\cal M}_{k}$ the set
\[
\CCs_{k}(c):= \biggl\{\mu\in\CCp\bigl([-a,a]\bigr)\dvtx  \int
_{-a}^{a}\Phi_{k}(x)\mu(dx)=c \biggr\}
\]
of all probability measures on the interval $[0,1]$ whose
moments up to the order $k$ coincide with $c =(c_1,\ldots, c_{k})^T$.
For $k = 2,3, \ldots$ and for a given point $(c_1,\ldots,
c_{k-1})^T \in{\cal M}_{k-1}$ we define $c^+_k = c^+_k(c_1,\ldots,
c_{k-1})$ and $c^-_k = c^-_k(c_1,\ldots,\break c_{k-1})$ as the largest and
smallest value of $c_k$ such that $(c_1,\ldots, c_k)^T \in\partial
{\cal M}_k$, that is,
\begin{eqnarray*}
c^{-}_k &=& \min\biggl\{ \int_{-a}^a x^k\mu(dx)\Big|\mu\in
{S}_{k-1}(c_1,\ldots, c_{k-1}) \biggr\},
\\
c^{+}_k &=& {\max} \biggl\{ \int_{-a}^a
x^k\mu(dx)\Big|\mu\in{S}_{k-1}(c_1,\ldots,
c_{k-1}) \biggr\}.
\end{eqnarray*}
Note that $c^-_k \le c_k \le c^+_k$ and that both inequalities are
strict if and only if $(c_1,\ldots, c_{k-1})^T \in\mathcal
{M}_{k-1}^0$ where $\mathcal{M}^0_{k-1}$ denotes the interior of the
set $\mathcal{M}_{k-1}$; see \citet{dettstud1997}.
For a moment point $c = (c_1,\ldots, c_n)^T$, such that $c = (c_1,\ldots,c_{n-1})^T$ is
in the interior of the moment space ${\cal M}_{n-1}$, the
canonical moments or canonical coordinates of the vector $c$ are
defined by $p_1=c_1$ and
%
\begin{equation}
\label{24} p_k = \frac{c_k - c^-_k}{c^+_k - c^-_k},\qquad k = 2,\ldots, n.
\end{equation}
Note that $p_k \in(0,1)$, $k=1,\ldots,n-1$ and $p_n\in\{0,1\}$
if and only if $(c_1,\ldots,\break c_{n-1}) \in\mathcal{M}^0_{n-1}$ and
$(c_1,\ldots, c_{n})^T \in\partial{\cal M}_n $.
In this case the canonical moments $p_i$ or order $i>n$ remain
undefined.\vadjust{\goodbreak}

We begin with a proof of the first part of Theorem~\ref{theo41}.
By Theorem~\ref{lem31} the determination of Bayesian $T$-optimal
discriminating designs can be obtained
by minimizing the linear optimality criterion
\[
\operatorname{tr} LM_{(2)}(\xi)
\]
for some appropriate matrix $L$, which is diagonal by the
symmetry of the prior distribution. A standard argument of optimal
design theory
shows that there exists a symmetric Bayesian $T$-optimal discriminating
design, say $\xi$,
for which the corresponding $2\times2$ matrix $M_{(2)}(\xi)$ is also
diagonal, that is,
\[
M_{(2)}(\xi)=\pmatrix{ a_{m-1}&0
\cr
0&a_m}.
\]
It now follows from Dette and Studden [(\citeyear{dettstud1997}),
Section 5.7], that for such a design, the elements in this matrix are
given by
%
\begin{equation}
\label{eq42} a_k (\xi) = (2a)^{2k} \prod
^k_{i=1} q_{2i-2}p_{2i-1}q_{2i-1}p_{2i},\qquad
k=m-1,m,
\end{equation}
where $q_0=1$, $q_i =1-p_i$ ($i\ge1$). Consequently, by Theorem
\ref{lem31} the Bayesian $T$-optimal discriminating design
problem is reduced to maximization of the function
%
\begin{equation}
\label{eq43}
\operatorname{tr} LM_{(2)}(\xi)=a_m(\xi) + \beta
a_{m-1}(\xi),
\end{equation}
where the quantities $a_m(\xi)$ are defined in
(\ref{eq42}), and $ \beta=\int b^2 \bar\pi(db)$ denotes the
second moment of the
prior distribution. This expression can now be directly
maximized in terms of the canonical moments, which gives
$p_{2m}=1$, $p_i=\frac{1}{2}$, $i=1,2,\ldots,2m-1$, $i\ne
2m-2$ and
\[
p_{2m-2}=\min\biggl\{ \frac{a^2+\beta}{2 a^2}, 1 \biggr\} =
\frac
{1+\beta_B}{2},
\]
where $\beta_B$ is defined in (\ref{betab}).
The corresponding design is uniquely determined and can be obtained
from Theorems 4.4.4 and 1.3.2 in \citet{dettstud1997}, which proves the
first part of the theorem.

For a proof of the second part we note that it follows from the proof
of Theorem~\ref{lem32} that the
standardized maximin $T$-optimal criterion
reduces to
%
\begin{equation}
\label{redopt} \inf_{b\in\B} \frac{a_m(\xi) +b^2 a_{m-1}(\xi) }{\bar
R(b)} \to\sup
_\xi,
\end{equation}
where $\bar R(b) $ is defined in (\ref{rquer}), that is,
\begin{eqnarray*}
\bar R (b)&=& \inf_{q_0,\ldots,q_{m-2}\in\mathbb{R}} \sup_{x \in
[-a,a]} \bigl|
x^m + b x^{m-1} + q_{m-2} x^{m-2} + \cdots+
q_1 x+q_0 \bigr|^2 \\
&=& \bar R(b,a).
\end{eqnarray*}
From (\ref{eq42}) it is obvious
that the canonical moments of a (symmetric) standardized maximin
$T$-optimal discriminating
design satisfy $p_{2m}=1$,
\[
p_i=\tfrac{1}{2},\qquad i=1,2,\ldots,2m-3,2m-1,
\]
and it remains to maximize (\ref{redopt}) with respect to the quantity
$p_{2m-2}$.\break A~straightforward calculation shows that the optimal value
of $p_{2m-2}$ is
determined by the condition $p_{2m-2}= 1-h^*$,
where $h^*$ is a solution of the problem
\[
\inf_{b\in\B} \frac{a^2h + b^2}{\bar R(b,a)a^2}(1-h)\to\max_{0\leq h\leq1/2}.
\]
The corresponding design is uniquely determined and can again be obtained
from Theorems 4.4.4 and 1.3.2 in \citet{dettstud1997}, which completes
the proof of Theorem~\ref{theo41}.
\end{pf}

\begin{remark}\label{remgaffke}
The structure of the Bayesian and standardized maximin $T$-optimal
designs determined in Theorem~\ref{theo41} is the same as the structure
of the $\phi _p$-optimal design for estimating the two coefficients
corresponding to the powers $x^m$ and $x^{m-1}$ in a polynomial
regression model of degree $m$ on the interval $[-a,a]$. More
precisely, it was shown in \citet{gaffke1987}, \citet{studden1989} (for
the interval $[-1,1]$) and in \citet{dettstud1997} (for arbitrary
symmetric intervals) that the designs minimizing
\[
\phi_p (\xi)= \bigl(\operatorname{tr}M_{(2)}^{-p}(\xi)
\bigr)^{1/p},\qquad -1<p\leq\infty,
\]
is given by the design $\xi_{m,\beta(p)}$ where $\beta(p)$ is the
unique solution of the equation
\[
\biggl( \frac{1-\beta}{2} \biggr)^{p+1} - a^{-2p} \beta=0
\]
in the interval $[0,1]$.
\end{remark}

\section{Some illustrative examples} \label{sec5}

In this section we illustrate the results in a few examples. We
restrict ourselves to the problem of discriminating between a constant
and the quadratic regression model on the interval $[-1,1]$.
Additionally, we construct robust designs for the situation considered
in Example~\ref{emax}. Further results for other models are
available from the authors.

\subsection{Standardized maximin $T$-optimal discriminating designs
for quadratic regression}

Consider the problem of discriminating between a constant and a
quadratic regression on the interval ${\cal X} = [-1,1]$.
As pointed out in \citet{brepinbra2005}, these models are of importance
for detecting dose response signals in phase II clinical trials.
If $\mathcal{B}=[-d,d]$, then it follows from Theorem \ref
{theo41} $(m=2)$ that
a standardized maximin $T$-optimal design is given by
\[
\xi^*_M=\pmatrix{ -1&0&1
\vspace*{2pt}\cr
\ds{\frac{1-h^*}{2}}& \ds{h^*}& \ds{
\frac{1-h^*}{2}}},
\]
where $h^*$ is a solution of the problem (\ref{eq41}). Due
to formula (3.9) in \citet{detmelshp2012} we have
%
\begin{equation}
\label{rbar} \bar R (b)= \bar{R}(b,1)=\cases{ \ds{ {1\over4}
\biggl(1+\frac{|b|}{2} \biggr)^4}, &\quad $|b|\leq2$,
\vspace*{2pt}\cr
b^2,&\quad $|b|\geq2$.}
\end{equation}
We define
\[
K(h,b)=\frac{h+b^2}{\overline R(b)} (1-h),
\]
and then the solution of the problem
\[
\max_{h \in[0,0.5] }\inf_{b \in[-d,d]} K(h,b)
\]
can be obtained by straightforward but tedious calculations, which are omitted
for the sake of brevity. For the solution one has
to distinguish three cases:
\begin{longlist}[(3)]
\item[(1)] If $0<d \leq\frac{1}{2}$
the minimum of the function $K(h,b)$ with respect to the variable $b$
is attained at the boundary of the interval $\B=[-d,d]$ and the
optimal value is given by $h^*=(1-d^2)/2$.
A typical situation is depicted in the left part of Figure~\ref{Pic01}.
%
\begin{figure}

\includegraphics{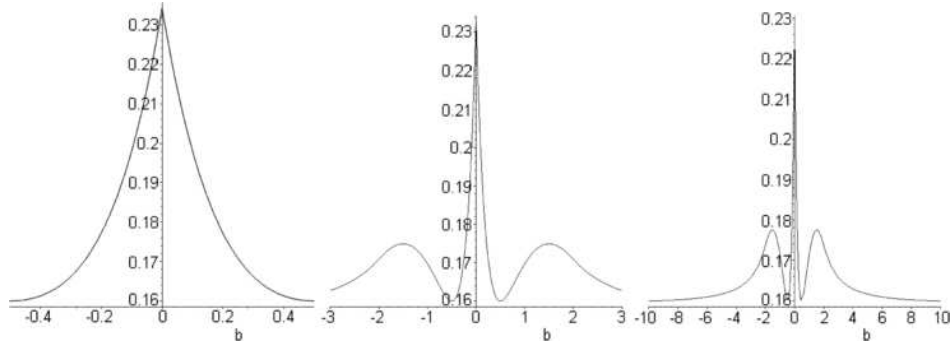}

\caption{The behavior of the function $K(h^*,b)$, for different
values of $d$. Left panel $d=1/2$, middle panel $d=2$, right panel
$d=10$.}\label{Pic01}
\end{figure}
A standardized maximin $T$-optimal discriminating design has masses
$ (1+d^2)/4$, $(1-d^2)/2 $ and $(1+d^2)/4$ at the points $-1$, $0$ and
$1$, respectively.
\item[(2)] In the case $1/2<d \leq5\sqrt{10}/4 $
the solution is given by
$h^*=3/8$, $b^*=1/2$. Therefore the design
with masses $5/16$, $3/8$ and $5/16$ at the points
$-1$, $0$ and $1$ is
a standardized maximin $T$-optimal discriminating design.
The behavior of the function $K(h^*,b)$ in this case is depicted in the
middle panel of Figure~\ref{Pic01}.
\item[(3)] In the case
$d\in[\frac{5\sqrt{10}}{4},\infty]$ the structure of the solution
changes again.
For this interval the optimal pair $h^*,b^*$ is obtained as a solution
of the system
\[
K(h,b)=K(h,d),\qquad \frac{\partial}{\partial b}K(h,b)=0,
\]
and we find by a direct calculation that
$ b^*$ is the unique root of the equation
\[
x^4+6x^3+\bigl(-2d^2+12
\bigr)x^2+\bigl(-16d^2+8\bigr)x+8d^2=0
\]
in the interval $[-4+2\sqrt{5},1/2]$. We have $h^*=
b^* -\frac{(b^*)^2}{2}$ and a standardized maximin $T$-optimal
discriminating design
has masses $1/2-b^*/2+(b^*)^2/4$, $b^*-(b^*)^2/2$, and $1/2-b^*/2+(b^*)^2/4$
at the points $-1$, $0$ and $1$, respectively. In the limiting case
$d=\infty$, that is, $\B= \mathbb{R}$, we have
$b^*=-4+2\sqrt{5}$, $h^*=-22+10\sqrt{5}$ and
a standardized maximin $T$-optimal discriminating design has masses
$23/2-5 \sqrt{5}$, $-22+10\sqrt{5}$, and $23/2-5 \sqrt{5}$
at the points $-1$, $0$ and~$1$, respectively.
A typical case for the function $K(h^*,b)$ in this case is depicted in the
right panel of Figure~\ref{Pic01} for $d=10$.
\end{longlist}

\subsection{Bayesian $T$-optimal discriminating designs for quadratic
regression}
\label{sec52}

For the Bayesian $T$-optimality criterion a prior has to be chosen,
and we propose to maximize an average of the efficiencies
\[
\int_{-a}^a \operatorname{eff}_T (\xi,b) \,db
\]
with respect to the uniform distribution on the interval $[-a,a]$. In
criterion (\ref{critvb}) this corresponds to an
absolute continuous prior with density proportional to
\[
f(b)=\cases{ \displaystyle \frac{3(2+a)^3}{16a(12+6a+a^2)}\frac{1}{\bar R(b)}, &\quad
$a\leq2$,
\cr
\displaystyle \frac{3a}{17a-6}\frac{1}{\bar R(b)}, &\quad $a\geq2$,}
\]
where $ \bar R(b) $ is defined in (\ref{rbar}) and the term depending
on $a$ is the corresponding
normalizing constant.
By direct calculations we obtain
\[
\int_{-a}^{a}b^2f(b)\,db=2\int
_{0}^{a}b^2f(b)\,db= \cases{
\displaystyle \frac{4a^2}{12+6a+a^2}, &\quad $a\leq2$,
\vspace*{2pt}\cr
\displaystyle \frac{6a^2-4a}{17a-6}, &\quad $a\geq2$.}
\]
In order to apply Theorem~\ref{theo41} we consider
$\beta_B=\min\{1,\int^a_{-a} b^2 f(b) \,db \} $, and again three
cases have to be considered:
\begin{longlist}[(3)]
\item[(1)]
If $0<a\leq2$ we have $ \beta_B= \frac{4a^2}{12+6a+a^2}$ and a
Bayesian
$T$-optimal discriminating design has masses $\frac
{5a^2+6a+12}{4(12+6a+a^2)}$, $\frac{-3a^2+6a+12}{2(12+6a+a^2)}$ and
$\frac{5a^2+6a+12}{4(12+6a+a^2)}$
at the points $-1$, $0$ and $1$.
\item[(2)] If $ 2\leq a\leq\frac{7+\sqrt{33}}{4} $
we have $ \beta_B= \frac{3a^2-4a}{17a-6}$, and a Bayesian
$T$-optimal discriminating design has masses $\frac
{6a^2+13a-6}{4(17a-6)}$, $\frac{-6a^2+21a-6}{2(17a-6)}$ and $\frac
{6a^2+13a-6}{4(17a-6)}$
at the points $-1$, $0$ and $1$.
\item[(3)] If $ \frac{7+\sqrt{33}}{4}\leq a $ we have $ \beta=1$
and a Bayesian
$T$-optimal discriminating design has masses $1/2$ and $ 1/2$
at the points $-1$ and $1$.
\end{longlist}

\subsection{Robust $T$-optimal discriminating designs for the
Michaelis--Menten and EMAX model}
\label{sec53}

\begin{figure}[t!]

\includegraphics{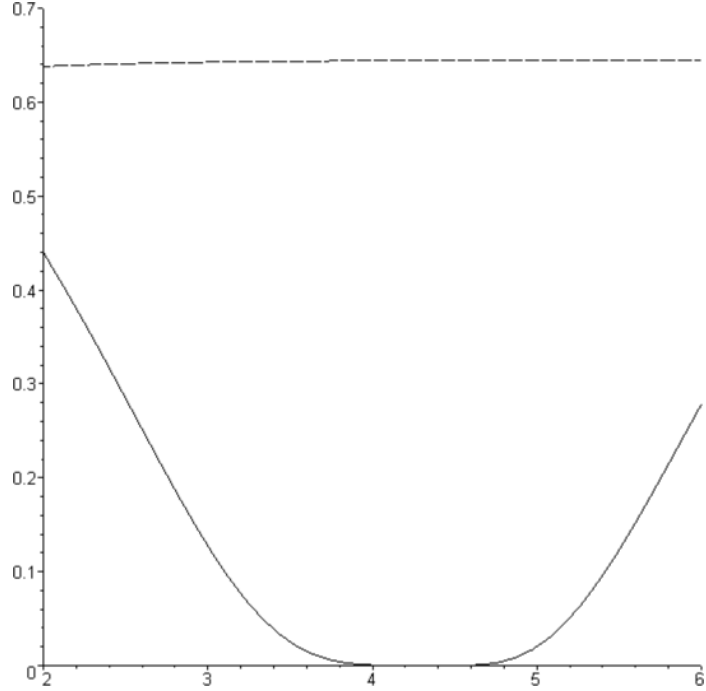}

\caption{$T$-efficiency (\protect\ref{teff}) of the standardized maximin
optimal discriminating design (dotted line) and the locally
$T$-optimal discriminating design for the Michaelis--Menten and EMAX
model (calculated under the assumption $\theta_{2,0}=-1/4$,
$\theta_{2,2}=1$, solid line). The ``true'' values are given by
$\theta_{2,0}=-1$, $\theta_{2,1}=1$ and the efficiencies depend on the
parameter $\theta_{2,2} \in(2,6)$.}\label{Pic5}
\end{figure}

In this section we briefly illustrate the application of the
methodology in the situation described in
Example~\ref{emax}, where the interest is in designs with good
properties for discriminating between the
Michaelis--Menten and EMAX model. We have calculated the standardized
maximin $T$-optimal
discriminating design for the Michaelis--Menten and EMAX model, where
the region for the
parameter $(\theta_{2,0},\theta_{2,1},\theta_{2,2})$ is given by
$[-1.1,-0.2] \times\{1\}
\times[2,6] $. The corresponding robust design is given by
\[
\xi= \pmatrix{ 1 & 1.36 & 2
\cr
0.410 & 0.205 & 0.385}.
\]
As pointed out in Example~\ref{emax}, the efficiency of locally
$T$-optimal discriminating designs can be low if some of the parameters
of the regression models have been misspecified, and in
Figure~\ref{Pic5} we compare the performance of the locally and
robust optimal discriminating designs if the true values are
$\theta_{2,0}=-1$, $\theta_{2,1}=1$ and
$\theta_{2,2} \in(2,6) $. We observe a substantial improvement
by the standardized maximin $T$-optimal discriminating design. Other scenarios
showed a similar picture and are not displayed for the sake of brevity.

%
%
%
\begin{figure}

\includegraphics{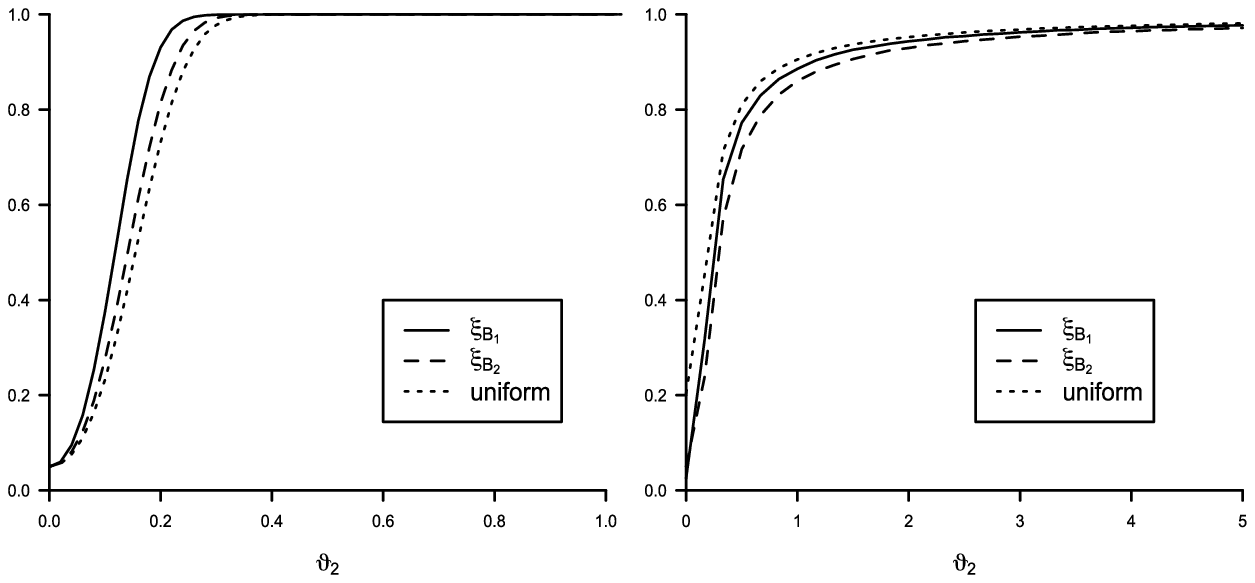}

\caption{The power of the $F$-test for the hypothesis (\protect\ref{sim2}).
Data are generated according to model (\protect\ref{simmod}). Left panel:
normal distributed errors; right panel: normal distributed errors
contaminated with $10\%$ Cauchy distributed errors.} \label{figsim}
\end{figure}
%

\subsection{Power and robustness}

In order to demonstrate the effect of the optimal design on the power
of the test for the corresponding hypothesis, we have conducted a small
simulation study comparing
Bayesian $T$-optimal discriminating designs and the commonly used
uniform designs with respect to their discrimination properties for the models
%
\begin{equation}
\label{sim1} \eta_1 (x,\theta_1) =
\theta_{1,0};\qquad \eta_2(x,\theta_2) =
\theta_{2,0} + \theta_{2,1} x + \theta_{2,2}
x^2,
\end{equation}
where the explanatory variable $x$ varies in the interval
$[-1,1]$.
For the construction of optimal discriminating designs we assume that
the ``true'' ratio of the coefficients of $x$ and $x^2$ is
an element of the interval $\mathcal{B}=[-1,1]$ or $\mathcal
{B}=[-3,3]$. It follows from Section~\ref{sec52} that the Bayesian $T$-optimal
discriminating designs with respect to the uniform distribution on the
interval $[-1,1]$ and $[-3,3]$ are obtained as
\[
\xi_{B_1} = \pmatrix{ -1 & 0 & 1
\vspace*{2pt}\cr
\frac{23}{76} &
\frac{15}{38} & \frac{23}{76} },\qquad \xi_{B_2} = \pmatrix{ -1 & 0
& 1
\vspace*{2pt}\cr
\frac{29}{60} & \frac{1}{30} & \frac{29}{60}},
\]
respectively. We assume that $60$ observations can be taken, which
yield to the ``realized'' designs:
\begin{itemize}
\item[$\bullet$] $\xi_{B1}$: $18$, $24$, $18$ observations at the
points $-1,0,1$ if $\mathcal{B}=[-1,1]$.
\item[$\bullet$] $\xi_{B_2}$: $29$, $2$, $29$ observations at the
points $-1,0,1$ if $\mathcal{B}=[-3,3]$.
\end{itemize}
For a comparison we use the uniform design:
\begin{itemize}
\item[$\bullet$] $6$ observations at the points $-1,-7/9,-5/9, \ldots,
5/9, 7/9, 1$.
\end{itemize}
In the left part of Figure~\ref{figsim} we show the simulated
rejection probabilities of the $F$-test for the hypothesis
%
\begin{equation}
\label{sim2} H_0\dvtx  \theta_{2,1} =\theta_{2,2} = 0
\end{equation}
(nominal level $5 \%$) in the model
%
\begin{equation}
\label{simmod} \eta_2(x,\theta_2) = 3 +
\tfrac{1}{2} \vartheta_2 x + \vartheta_2x^2,
\end{equation}
where the errors are centered normal distributed with variance $\sigma
^2=0.5$ (note that this means
that the ``true'' ratio of the coefficients of $x$ and $x^2$ is given
by $b=\frac{\theta_{2,1}}{\theta_{2,2}}=1/2$). All results are based
on 250,000 simulation runs. We observe a notable
improvement with respect
to the power of the $F$-test if the experiments are conducted according
to the
$T$-optimal discriminating designs. The Bayesian $T$-optimal discriminating
design with respect to the uniform distribution on $\mathcal
{B}=[-1,1]$ yields a larger power than the Bayesian optimal design with
respect to the uniform distribution on
the interval $[-3,3]$. This corresponds to intuition because this
design uses more precise and correct information regarding the unknown
ratio $b=\theta_{2,1}/\theta_{2,2}$. In fact, using Bayesian $T$-optimal
discriminating designs with respect to smaller intervals containing the
``true'' value $b=1/2$ yields even more powerful tests (these results
are not depicted for the sake of brevity).

It was pointed out by a referee that it might be of interest to
investigate the sensitivity of the $F$-test for the different designs
with respect to influential observations. For this purpose we have
performed the same simulation where $10 \%$ of the normal distributed
errors are replaced by Cauchy distributed random variables. The
corresponding results are shown in the right part of Figure~\ref{figsim}, and
the results change substantially. We observe a loss in power for all
three designs. Under the null hypothesis the Bayesian $T$-optimal discriminating
designs yield a slightly conservative test while the $F$-test based on
the uniform design rejects the null hypothesis too often. Because of
continuity of the power function this phenomenon is also observed for
other values of $\vartheta_{2}$. On the other hand, for large values
of $\vartheta_{2}$ the Bayesian $T$-optimal discriminating design
$\xi_{B_1} $ with respect to the uniform distribution on the interval
$[-1,1]$ and the uniform design yield a similar power of the $F$-test, while
a slightly lower power is observed for the $F$-test based on the
Bayesian $T$-optimal discriminating design $\xi_{B_2} $ with respect
to the uniform distribution on the interval $[-3,3]$.

\begin{figure}

\includegraphics{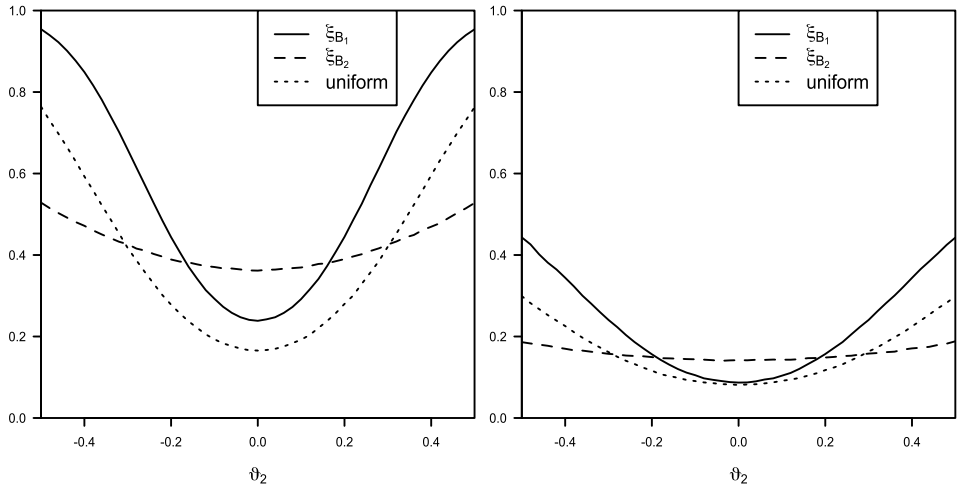}

\caption{The power of the $F$-test for the hypothesis (\protect\ref{sim2}).
Data are generated according to model (\protect\ref{simmod1}). Left panel:
normal distributed errors; right panel: normal distributed errors
contaminated with $10\%$ Cauchy distributed errors.}
\label{figsim1}
\end{figure}

Model (\ref{simmod}) keeps the ratio of the coefficients of $x$ and
$x^2$ constant and in the second example
of this section we consider
an alternative data generating model, that is,
%
\begin{equation}
\label{simmod1} \eta_2(x,\theta_2) = 3 +
\tfrac{1}{8} x + \vartheta_2 x^2,
\end{equation}
where $\vartheta_2$ varies in the interval $[-0.5,0.5]$, which means
that the ``true'' ratio $b$ of
the coefficients of $x$ and $x^2$ varies in $\real\setminus
[-0.25,0.25]$. In particular the
ratio of the coefficients of $x$ and $x^2$
in model (\ref{simmod1}) can attain values which are not contained in
the set $\mathcal{B}$ used for the
construction of the Bayesian $T$-optimal discriminating designs.
Again $60$ observations are generated according to the designs
specified in the previous paragraph and the corresponding results are
depicted in Figure~\ref{figsim1}. Note that none
of the values $ \vartheta_2 \in[-0.5,0.5]$ corresponds to the null
hypothesis and consequently the curves show only
the power under certain alternatives. We observe from the left panel in
Figure~\ref{figsim1} that for normally distributed
%
errors the Bayesian $T$-optimal discriminating design $\xi_{B_1} $
with respect to the uniform distribution on the interval $[-1,1]$ yields
uniformly more power than the uniform design. On the other hand the
Bayesian $T$-optimal discriminating design $\xi_{B_2} $ with respect
to the uniform distribution on the interval $[-3,3]$ is preferable to
the uniform
design if $\vartheta_2 \leq0.3$, while for larger values of
$\vartheta_2$ the $F$-test based on the uniform design is more powerful.
Moreover, if
$\vartheta_2 \leq0.1$ the design $\xi_{B_2} $ is even better than
the design $\xi_{B_1} $. This corresponds to intuition, because
the design $\xi_{B_2} $ is very close to the optimal design for
discriminating between a constant and a linear regression model,
which puts equal masses at the points $-1$ and $1$. For Cauchy
distributed errors we observe a very similar behavior, where
the power is smaller due to the contamination of the normal distribution.

\section*{Acknowledgments}

The authors thank two unknown referees for their constructive comments
on an earlier version of this manuscript and Martina Stein, who typed
parts of this manuscript with considerable technical expertise. We are
also grateful to Katrin Kettelhake for numerical assistance. This paper
was initiated at the Isaac Newton Institute for Mathematical Sciences
in Cambridge, England, during the 2011 programme on the Design and
Analysis of Experiments.



\printaddresses

\end{document}